\documentclass[12pt,twoside,reqno]{amsart}

\usepackage{bm}

\usepackage{amssymb,amsmath,amstext,amsthm,amsfonts,amscd,xcolor}

\usepackage{dsfont}

\usepackage{amsthm, enumerate}

\usepackage[ansinew]{inputenc}
\usepackage{graphicx}
\usepackage[mathscr]{eucal}
\usepackage{hyperref}

\newcommand{\R}{\mathbb{R}}

\newcommand{\N}{\mathbb{N}}
\newcommand{\Z}{\mathbb{Z}}

\newcommand{\T}{\mathbb{T}}
\newcommand{\SL}{{\rm SL}}

\newcommand{\Pp}{\mathbb{P}}

\newcommand{\EE}{\mathbb{E}}
\newcommand{\Escr}{\mathscr{E}}

\newcommand{\Sigmah}{{\Omega}}

\theoremstyle{plain}
\newtheorem{theorem}{Theorem}[section]
\newtheorem{proposition}{Proposition}[section]

\theoremstyle{definition}

\newtheorem{definition}{Definition}[section]
\newtheorem*{theorem*}{Theorem}

\theoremstyle{definition}

\numberwithin{equation}{section}

\newcommand{\abs}[1]{\left| #1 \right|} % absolute value
\newcommand{\norm}[1]{\left\|#1\right\|} % norm
\newcommand{\normtwo}[1]{% Peter Grill norm @tex.stackexchange.com
{\left\vert\kern-0.25ex\left\vert\kern-0.25ex\left\vert #1
    \right\vert\kern-0.25ex\right\vert\kern-0.25ex\right\vert} }

\newcommand{\ep}{\epsilon}

\newcommand{\om}{\omega}

\makeatletter
\newsavebox\myboxA
\newsavebox\myboxB
\newlength\mylenA

\newcommand*\xoverline[2][0.75]{%
    \sbox{\myboxA}{$\m@th#2$}%
    \setbox\myboxB\null% Phantom box
    \ht\myboxB=\ht\myboxA%
    \dp\myboxB=\dp\myboxA%
    \wd\myboxB=#1\wd\myboxA% Scale phantom
    \sbox\myboxB{$\m@th\overline{\copy\myboxB}$}%  Overlined phantom
    \setlength\mylenA{\the\wd\myboxA}%   calc width diff
    \addtolength\mylenA{-\the\wd\myboxB}%
    \ifdim\wd\myboxB<\wd\myboxA%
       \rlap{\hskip 0.5\mylenA\usebox\myboxB}{\usebox\myboxA}%
    \else
        \hskip -0.5\mylenA\rlap{\usebox\myboxA}{\hskip 0.5\mylenA\usebox\myboxB}%
    \fi}
\makeatother

\newcommand{\Proj}{\mathbb{P}(\R^m)}

\newcommand{\dist}{{\rm dist}}

\newcommand\restr[2]{{% we make the whole thing an ordinary symbol
  \left.\kern-\nulldelimiterspace % automatically resize the bar with \right
  #1 % the function
  \vphantom{\big|} % pretend it's a little taller at normal size
  \right|_{#2} % this is the delimiter
  }}
\newcommand{\Prob}{\mathrm{Prob}}
\newcommand{\Ascr}{\mathscr{A}}
\newcommand{\Gscr}{\mathscr{G}}

\newcommand{\supp}{\mathrm{supp}}

\newcommand{\Qop}{\mathcal{Q}}

\newcommand{\alfa}{\mathfrak{a}}
\newcommand{\muh}{{\nu}}
\newcommand{\nuh}{{\tilde \nu}}

\newcommand{\Xplus}{{X^+}}
\newcommand{\fle}{L_1}

\newcommand{\Bcal}{\mathcal{B}}

\title{Randomness versus quasi-periodicity}

\date{}

\begin{document}

\author[A. Cai]{Ao Cai}
\address{School of Mathematical Sciences\\
Soochow University\\
China
}
\email{acai@suda.edu.cn}

\begin{abstract}
This paper serves as an extended road map for our long-term project ``Mixed Random-quasiperiodic Cocycles''\cite{CDK-paper2,CDK-paper1,CDK-paper3,CDK-paper4,CDK-paper5} with Pedro Duarte and Silvius Klein. Despite exhibiting totally different natures, the random world and the quasi-periodic one may still have potential relations that we are keen to reveal. This was inspired by Jiangong You's intriguing question on the stability of the Lyapunov exponent of quasi-periodic Schr\"odinger operators under random noise in 2018.
\end{abstract}

\maketitle

\tableofcontents

\section{Introduction}\label{intro}
The story began during an international conference on dynamical systems held in Nanjing University in 2018 when I met Duarte for the first time and Klein for the second time who both became my postdoctor supervisors afterwards . Back then, I just finished my PhD with my dissertation ``Reducibility of finitely differentiable quasi-periodic cocycles and its applications'' based on Kolmogorov-Arnold-Moser (KAM) theory under the supervision of You. Obviously, I was a purely quasi-periodic person~\cite{CaiAC,CCYZ,CaiGe,CaiWang} while Duarte and Klein had already many collaborations on random cocycles as well as quasi-periodic ones, see the two excellent books~\cite{DKLEbook,DK-31CBM} and the references therein. 

Definitely, You's beautiful question motivated the two parallel rays to intersect and resonate intensively.

Digressions aside, You's original question states:``what is the behavior of the Lyapunov exponent of the quasi-periodic Schr\"odinger operator (cocycle) perturbed by an i.i.d. random noise?'' 

More precisely, consider the one dimensional discrete Schr\"odinger operator defined on $\ell^2(\Z)$:
$$
(H_{W,V,\alpha,\theta}u)_n=u_{n+1}+u_{n-1}+[V(\theta+n\alpha)+W_n]u_n
$$
where $V$ is a real-valued potential function and $\{W_n\}_{n\in \Z}$ is an i.i.d. sequence driven by some probability measure $\mu$ on $\R$. If we put a coupling $\epsilon>0$ before $W_n$ and let $\epsilon$ go to $0$, what will be the behavior of the Lyapunov exponent as a function of $\epsilon$ and $\mu$? Namely, You is concerned with the stability of the Lyapunov exponent of quasi-periodic Schr\"odinger cocycles under random perturbation when the magnitude of the randomness vanishes.

This great question motivates us to use everything we have as well as develop various new methods, trying to build the bridge between the quasi-periodicity and the randomness. For the time being, we still can not give a complete answer about it but partial results are obtained. In fact, along the whole way of attempting to answer it, we already acquired many fruitful results.

In this short paper, we provide five sections of results contained respectively in~\cite{CDK-paper2,CDK-paper1,CDK-paper3,CDK-paper4,CDK-paper5} where relations and comparisons of our theorems with the existing literature can be found. Since our main purpose is to convey ideas and state results with no proof at all, we decide to leave the statements of main theorems in their corresponding sections in order to avoid repetition.

This paper is organized as follows. In section 2, we introduce the concept of mixed random-quasiperiodic cocycles. In section 3, we state our Furstenberg positivity criterion of the maximal Lyapunov exponent of mixed models and give some applications in Mathematical Physics. Section 4 is an interlude which mainly contains an abstract large deviations type theorem. In section 5, we talk about the H\"older continuity of the Lyapunov exponent. Finally, section 6 is devoted to the stability of the Lyapunov exponent.

\section{Mixed random-quasiperiodic cocycles}\label{generalities}
As preliminaries, in this section we introduce the concept of mixed random-quasiperiodic cocycles. While the basic definitions could not be omitted much, all the details of properties will be intentionally excluded. Interested readers are kindly invited to~\cite{CDK-paper1} for precise statements and proofs.

\subsection*{The base dynamics}  Let $(\Sigmah, \Bcal)$ be a standard Borel space
and let $\muh \in \Prob_c (\Sigmah)$ be a compactly supported Borel probability measure on $\Sigmah$. Regarding $(\Sigmah, \muh)$ as a space of symbols, we consider the corresponding (invertible) Bernoulli system $\left(X, \sigma, \muh^\Z \right)$, where $X := \Sigmah^\Z$ and $\sigma \colon X \to X$ is the (invertible) Bernoulli shift: for $\om = \{ \om_n \}_{n\in\Z}  \in X$,  
$\sigma \om  :=  \{ \om_{n+1} \}_{n\in\Z}$. Consider also its non invertible factor on $X^+ := \Sigmah^\N$.

Let $\T^d = \left(\R/\Z \right)^d$ be the torus of dimension $d$, and denote by $m$ the Haar measure on its Borel $\sigma$-algebra.

Given a continuous function $\alfa \colon \Sigmah \to \T^d$, the skew-product map
\begin{equation}\label{base map}
	f \colon X \times \T^d \to X  \times \T^d \, , \quad f (\om, \theta) := \left( \sigma \om, \theta + \alfa ( \om_0) \right)
\end{equation}
will be referred to as a mixed random-quasiperiodic (base) dynamics. 

This map preserves the measure $\muh^\Z\times m$ and it is the natural extension of the non-invertible map on
$X^+ \times \T^d$ which preserves the measure $\muh^\N\times m$ and is defined by the same expression.

We call the measure $\muh$ ergodic, or ergodic with respect to $f$ when the mixed random-quasiperiodic system $\left( X \times \T^d, f, \muh^\Z\times m \right)$ is ergodic. See~\cite[Section 2]{CDK-paper1} for various characterizations of the ergodicity. In the same paper, a uniform convergence of Birkhoff sums of continuous observables to their space averages was established~\cite[Lemma 2.5]{CDK-paper1}. This was further used, along with a stopping time argument, to prove a uniform base LDT theorem for continuous observables depending on finitely many coordinates~\cite[Theorem 2.4]{CDK-paper1}.

For our interest, before introducing the fiber dynamics, we are going to specify the $\Omega$ as follows.

\subsection*{The group of quasiperiodic cocycles} Given a  frequency $\alpha \in \T^d$, let $\tau_\alpha (\theta) = \theta + \alpha$ be the corresponding ergodic translation on $\T^d$. Consider $A \in C^0 (\T^d, \SL_m (\R))$ a continuous matrix valued function on the torus. A quasiperiodic cocycle is a skew-product map of the form 
$$\T^d \times \R^m \ni (\theta, v) \mapsto \left( \tau_\alpha (\theta),  A(\theta) v \right) \in \T^d \times \R^m \, .$$

This cocycle can thus be identified with the pair $(\alpha, A)$. Consider the set
$$\Gscr=\Gscr(d,m):= \T^d\times C^0(\T^d,\SL_m(\R))$$
of all quasiperiodic cocycles.  

This set is a Polish metric space when equipped with the product metric (in the second component we consider the uniform distance). The space $\Gscr$ is also a group, and in fact a topological group, with the natural composition and inversion operations
\begin{align*}
	(\alpha,A)\circ (\beta, B) &:= (\alpha+\beta, (A\circ \tau_\beta) \, B)  \\
	(\alpha,A)^{-1} &:= ( -\alpha, (A\circ \tau_{-\alpha})^{-1} ) \, .
\end{align*}

Given $\muh \in \Prob_c (\Gscr)$ let $\om = \left\{ \om_n \right\}_{n\in\Z}$, $\om_n = (\alpha_n, A_n)$ be an i.i.d. sequence of random variables  in $\Gscr$  with law $\muh$. Consider the corresponding multiplicative process in the group $\Gscr$
\begin{align*}
	\Pi_n & = \om_{n-1} \circ \cdots \circ \om_1 \circ \om_0 \\
	& = \left( \alpha_{n-1} + \cdots + \alpha_1 + \alpha_0, \, ( A_{n-1} \circ \tau_{\alpha_{n-2} + \cdots + \alpha_0}) \cdots (A_1 \circ \tau_{\alpha_0} ) \, A_0   \right) \, .
\end{align*}

In order to study this process in the framework of ergodic theory, we model it by the iterates of a  linear cocycle.

\subsection*{The fiber dynamics} Given $\muh \in \Prob_c (\Gscr)$, let $\Sigmah \subset \Gscr$ be a closed subset such that $\Sigmah \supset \supp \muh$. 
We regard $(\Sigmah, \muh)$ as a space of symbols and consider the shift $\sigma$ on the space $X := \Sigmah^\Z$ of sequences $\om = \left\{ \om_n \right\}_{n\in\Z}$ endowed with the product measure $\muh^Z$ and the product topology (which is metrizable). 
The standard projections 
\begin{align*}
	\alfa \colon \Sigmah \to \T^d, & \qquad \alfa (\alpha, A) = \alpha \\
	\Ascr \colon \Sigmah \to C^0 (\T^d, \SL_m (\R)), & \quad \quad \Ascr (\alpha, A) = A
\end{align*}
determine the linear cocycle $F=F_{(\alfa,\Ascr)} \colon X\times\T^d\times\R^m \to X\times\T^d\times\R^m$ defined by
$$ F(\omega,\theta, v) := \left(\sigma \omega, \theta+\alfa(\omega_0), \Ascr(\omega_0)(\theta)\, v \right)  .$$
The non-invertible version of this map is defined similarly on $\Xplus\times \T^d\times \R^m$, where $\Xplus=\Sigmah^{\N}$.

Thus the base dynamics of the cocycle $F$ is the mixed random-quasiperiodic map $f$ defined above,
$$ X \times \T^d \ni (\om, \theta) \mapsto f (\om, \theta) := \left(\sigma \om, \theta + \alfa (\om_0) \right) \in X \times \T^d ,$$
while the fiber action is induced by the map
$$X \times \T^d \ni ( \om, \theta ) \mapsto \Ascr (\om, \theta) := \Ascr (\om_0) (\theta) \in \SL_m (\R) .$$

The skew-product $F$ will then be referred to as a {\em mixed random-quasiperiodic cocycle}. The space of mixed cocycles $F = F_{(\alfa, \Ascr)}$ is a metric space with the uniform distance
$$
\dist \left( (\alfa,\Ascr), \, (\alfa',\Ascr') \right)= \norm{\alfa-\alfa'}_0 + \norm{\Ascr-\Ascr'}_0 \, .
$$

For $\om = \{ \om_n \}_{n\in\Z} \in X$ and $j \in \N$ consider the composition of random translations
\begin{align*}
	\tau_\om^j & := \tau_{\alfa(\om_{j-1})} \circ \cdots \circ \tau_{\alfa(\om_0)} 
	= \tau_{ \alfa(\om_{j-1}) + \cdots +  \alfa(\om_0) } = \tau_{\alfa (   \om_{j-1} \circ \cdots \circ \, \om_0  ) } \, .
\end{align*}

The iterates of the cocycle $F$ are then given by
$$F^n (\om, \theta, v) = \left( \sigma^n \om, \tau_\om^n (\theta), \, \Ascr^n (\om) (\theta) v \right) ,$$
where
\begin{align*}
	\Ascr^n (\om) & = \Ascr \left( \om_{n-1} \circ \cdots \circ \om_1 \circ \om_0 \right) \\
	& = \left( \Ascr (\om_{n-1}) \circ \tau_\om^{n-2}  \right) \, \cdots \,  \left( \Ascr (\om_{1}) \circ \tau_\om^{0}  \right) \, \Ascr (\om_0) \, .
\end{align*}
Thus  $\Ascr^n(\omega)$ can be interpreted as a random product of
quasiperiodic cocycles, driven by the measure $\muh$ on the group $\Gscr$ of such cocycles. For convenience we also denote $\Ascr^n (\om, \theta) := \Ascr^n (\om) (\theta)$.

\medskip

By Kingman's subadditive ergodic theorem, the limit of the sequence $\displaystyle \frac{1}{n} \, \log \norm{ \Ascr^n (\om) (\theta)}$ as $n\to\infty$ exists for $\muh^\Z \times m$ a.e. $(\om, \theta) \in X \times \T^d$. Recall that for simplicity we call the measure $\muh$ ergodic if the base dynamics $f$ is ergodic w.r.t. $\muh^\Z \times m$. In this case, the limit is a constant depending only on the measure $\muh$ and it is called the first (or maximal) Lyapunov exponent of the cocycle $F$, which we denote by $L_1 (\muh)$. Thus
\begin{align*}
	L_1 (\muh) & = \lim_{n\to\infty} \frac{1}{n} \, \log \norm{ \Ascr^n (\om) (\theta)} \quad \text{for} \quad \muh^\Z \times m \text{ a.e. } (\om, \theta) \\
	& = \lim_{n\to\infty} \, \int_{X \times \T^d}  \frac{1}{n} \, \log \norm{ \Ascr^n (\om) (\theta)} \, d ( \muh^\Z \times m ) \, .
\end{align*}
The other Lyapunov exponents are similarly defined via the corresponding singular values.

Using the available uniform base LDT theorem as well as the subadditive ergodic theorem, we proved a uniform fiber upper LDT estimates~\cite[Theorem 3.1]{CDK-paper1} which readily implies the upper semi-continuity of the Lyapunov exponent with respect to the Wasserstein distance. Moreover, at this level of generality, the uniform fiber lower LDT can not necessarily hold~\cite[Remark 3.2]{CDK-paper1}. Otherwise we get continuity for free but there are counter examples cooked up based on Wang-You~\cite{WangYou}.

After the basic concepts and properties are introduced, we may start building various results upon them.

\section{Positivity of the Lyapunov exponent}\label{generic}

\subsection*{Furstenberg positivity criterion}
The first interesting result is the following criterion for the positivity of the maximal Lyapunov exponent of mixed random-quasiperiodic cocycles.

Before stating the theorem, we introduce two useful definitions.

\begin{definition}[See Definition 3.16 in~\cite{FurmanSurvey}]
	\label{non-compact/strongly irreducible, Zariski: Gscr}
	
	Let $\Gscr_0\subset \Gscr$ be a closed subgroup.
	We say that $\Gscr_0$ is respectively \textit{non-compact},
	\textit{strongly irreducible} or \textit{Zariski dense} when  there exists no measurable function $M\colon \T^d\to \SL_m(\R)$ such that
	$C:\Gscr_0\times\T^d\to \SL_m(\R)$,
	$$  C((\beta, B), \theta):= M(\theta+\beta)^{-1}\, B(\theta) \,M(\theta) $$
	for any $(\beta, B)\in \Gscr_0$ and $m$-a.e. $\theta$, takes values in a proper closed subgroup  $G_0\subseteq \SL_m(\R)$ that is compact,   virtually reducible or  algebraic, respectively.

	Given measure $\muh\in\Prob_c(\Gscr)$ we denote by $\Gscr_\muh$ the closed subgroup of $\Gscr$ generated by $\supp(\muh)$.
	This measure $\muh$ is called  \textit{non-compact},
	\textit{strongly irreducible} or \textit{Zariski dense} when
	the closed subgroup $\Gscr_\mu$ satisfies the same properties.
\end{definition}

\begin{theorem}\label{thm intro positivity}
	Let $\muh\in \Prob_c(\Gscr)$ and assume that $\muh$ is ergodic, non-compact and strongly irreducible.
	Then $\fle(\muh)>0$.
\end{theorem}
 This naturally extends the classical Furstenberg's positivity criterion for products of random matrices \cite{Fur} and its proof is indeed in the spirit of Furstenberg's original one. Moreover, it also has intriguing applications in Mathematical Physics.
 
 \subsection*{Applications to Schr\"odinger operators} Consider the Schr\"odinger operator $H (\theta)$ on $l^2 (\Z)$ defined by
 \begin{equation}\label{op1}
 	(H (\theta)  \psi)_n = - \psi_{n+1} - \psi_{n-1} + \left( v (\theta + n \alpha) + w_n \right) \, \psi_n \quad \forall n \in \Z ,
 \end{equation}
 for all $\psi = \{\psi_n\}_{n\in\Z} \in l^2 (\Z)$.
 
 The i.i.d. sequence $\{w_n\}$ can be interpreted as a random perturbation of the quasiperiodic potential $\{v_n (\theta)\}$, where $v_n (\theta) := v (\theta + n \alpha)$.
 
 \medskip
 
 We may instead randomize the frequency, that is, given an i.i.d sequence $\{\alpha_n\}_{n\in\Z}$ of random variables with values on the torus $\T^d$, 
 consider the discrete Schr\"odinger operator
 \begin{equation}\label{op2}
 	(H (\theta)  \psi)_n = - \psi_{n+1} - \psi_{n-1} + v (\theta + \alpha_0 + \cdots + \alpha_{n-1} )  \, \psi_n \quad \forall n \in \Z .
 \end{equation}

 We may certainly both randomize the frequency and randomly perturb the resulting potential (independently from each other or not), thus obtaining the mixed Schr\"odinger operator
 \begin{equation}\label{op3}
 	(H (\theta)  \psi)_n = - \psi_{n+1} - \psi_{n-1} + \left( v_n (\theta) + w_n \right) \, \psi_n \quad \forall n \in \Z,
 \end{equation}
 where
 $$v_n (\theta) = v (\theta + \alpha_0 + \cdots + \alpha_{n-1} ) \, .$$
 
 For these three models, we have respectively three applications.
 
\begin{proposition}\label{intro prop1}
	Consider the Schr\"odinger operator~\eqref{op1} with randomly perturbed quasiperiodic potential and the corresponding cocycles driven by the measures
	$\muh_E := \delta_\alpha \times \int_\R \delta_{P (\om) S_E} \, d \rho (\om)$
	where $\rho \in \Prob_c (\R)$.
	
	If $\supp (\rho)$ has more than one element, then $L_1 (\muh_E) > 0$ $\forall E \in \R$. 
\end{proposition}

We now consider the case of randomly chosen frequencies.

\begin{proposition}\label{intro prop2}
	Consider the Schr\"odinger operator~\eqref{op2} with randomly perturbed frequencies and the corresponding cocycles driven by the measures $\muh_E :=  \mu   \times  \delta_{S_E} $
	where $\mu \in \Prob (\T^d)$.
	
	If there are two frequencies $\alpha, \beta \in \supp (\mu)$ such that $\beta - \alpha$ is an ergodic translation on $\T^d$, and if the potential function $v (\theta)$ is analytic and non-constant (or, more generally, if the continuous functions $v (\theta)$ and $v (\theta + \beta - \alpha)$ are transversal), then $L_1 (\muh_E) > 0$ $\forall E \in \R$. 
\end{proposition}

\begin{proposition}\label{intro prop3}
	Consider the Schr\"odinger operator~\eqref{op3} with randomly perturbed frequencies and quasiperiodic potential and the corresponding cocycles driven by the measures
	$\muh_E= \int_{\T^d \times \R} \delta_{(\alpha, \, P (\om) S_E)} \, d \eta (\alpha, \om)$
	where $\eta \in \Prob_c (\T^d \times \R)$.
	
	If $\supp (\eta)$ contains two points $(\alpha, \om_1)$ and $(\alpha, \om_2)$ with $\alpha$ rationally independent and $\om_1 \neq \om_2$, then $L_1 (\muh_E) > 0$ $\forall E \in \R$. 
\end{proposition}

It is not hard to see that randomness dominates quasi-periodicity in the sense of positive Lyapunov exponents. For example,  for any quasi-periodic Schr\"odinger operator in the almost reducible regime (small coupling constant, Diophantine frequency and smooth enough), it has purely absolutely continuous spectrum with zero Lyapunov exponents~\cite{CaiAC}\cite{CCYZ}. However, when the randomness comes in, the new mixed cocycle has positive Lyapunov exponent which is consistent with the random case. Indeed, by Kotani theory, the spectral type also changes completely from absolutely continuous spectrum to singular spectrum.

Note that the same paper also provides an average uniform convergence to the top Lyapunov exponents in the Oseledets theorem, which is an essential preparation for proving the regularity, namely H\"older continuity, of the Lyapunov exponent under generic assumptions.

Although it would be quite natural to continue with a new section on regularity, let us have a break here, which will proved to be worthy in a moment.

\section{Interlude: an abstract LDT theorem }\label{interlude}
As is well known, one very powerful tool to proving the continuity of the Lyapunov exponent is the so-called large deviations type (LDT) theorem. In fact, Duarte and Klein~\cite{DKLEbook} built an abstract continuity theorem (ACT) assuming uniform fiber LDT holds. Therefore, it is wise and natural to focus on LDT before going further.

In this section, we introduce an abstract LDT theorem for strongly mixing Markov systems obtained in~\cite{CDK-paper3}. Let us start by several basic definitions.

\subsection*{Markov systems}
\begin{definition}
	A stochastic dynamical system (SDS) is any continuous map
	$K \colon M\to \Prob(M)$ on a Polish metric space $M$.
\end{definition}

Given any $x\in M$ and $E\subset M$ measurable, $K_x(E)$ represents the probability of going from $x$ into $E$. The iterates are
$$
K^{n+1}_x(E)=\int_M K^n_y(E)dK_x(y).
$$

\begin{definition}
	A probability measure $\nu \in \Prob (M)$ is called a $K$-stationary if $$\nu = K\ast \nu := \int K_x\, d\nu(x).$$
\end{definition}

\begin{definition}
The Markov operator induced by an SDS $K$:
$\Qop=\Qop_K \colon L^\infty(M,d\nu)\to L^\infty(M,d\nu)$ is defined by 
	$$ (\Qop \varphi)(x):= \int_M \varphi(y)\, dK_x(y) .$$
\end{definition}

$K$ is also called the Markov kernel. Similarly,
$$
(Q^n\varphi)(x)=\int_M \varphi(y)dK_x^n(y)
$$

\begin{definition}\label{marsys}
	A Markov system is a tuple $(M,K,\nu,\Escr)$  where
	\begin{enumerate}
		\item $M$ is a Polish metric space,
		
		\item $K \colon M\to \Prob(M)$ is an SDS,
		
		\item  $\nu\in\Prob(M)$ is a $K$-stationary measure,
		
		\item  $\Escr=\left(\Escr,\norm{\cdot}_\Escr\right)$ is a Banach subspace
		of $L^\infty(M,d\nu)$ such that the inclusion $\Escr\subset L^\infty(M,d\nu)$ and the action of $\Qop$ on $\Escr$ are both continuous. In other words there are constants $M_1<\infty$ and $M_2<\infty$ such that
		$\norm{\varphi}_\infty\leq M_1\, \norm{\varphi}_\Escr$ and
		$\norm{\Qop\varphi}_\Escr\leq M_2\, \norm{\varphi}_\Escr$,
		for all $\varphi\in \Escr$. 
	\end{enumerate}
\end{definition}

\subsection*{Strong mixing}

Now let us recall by order three definitions of strong mixing.

Firstly, the strongest version of mixing is:
\begin{definition}[Meyn-Tweedie~\cite{Markov-chains-book}]
	The Markov system is called uniformly ergodic if for any $x\in M$
	$$
	K_x^n \to \nu \quad \text{as}\quad n\to \infty.
	$$
\end{definition}

This condition is equivalent to the following: there exist some $C<\infty$ and $\sigma\in (0,1)$ such that 
	$$
	\norm{\Qop^n \varphi -\int_M \varphi \, d\nu}_\infty \leq C\sigma^n \norm{\varphi}_\infty
	$$
for all $n\in \N$ and $\varphi \in L^\infty(M,d\nu)$.

Now, let us state a weaker version.
\begin{definition}[Duarte-Klein\cite{DKLEbook}]
	There are constants $C < \infty$ and $\sigma \in (0, 1)$ such that
		$$ \norm{\Qop^n \varphi-\int_M \varphi\, d\nu}_\Escr \leq C \,\sigma^n \norm{\varphi}_\Escr  $$
	for all $\varphi \in \Escr$ and $n \in \N$.
\end{definition}

Note that the two definitions above lead to exponential LDT.

Finally, we may introduce the weakest version.
\begin{definition}
	A Markov system $(M, K,\nu,\Escr)$ is  called strongly mixing with power mixing rate  if
	there are constants  $C<\infty$ and $p > 0$ such that for all $\varphi\in\Escr$ and  $n\in\N$,
	
		\begin{equation*}
			\norm{\Qop^n \varphi-\int_M \varphi\, d\nu}_\infty \leq C \, \frac{1}{n^p} \, \norm{\varphi}_\Escr  \, . 
	\end{equation*}
\end{definition}

\subsection*{Large deviations type theorem}
As we know, the LDT for a measure-preserving dynamical system (MPDS) can be seen as a measure version of Birkhoff ergodic theorem, so the Birkhoff sum is involved:
	for all $\epsilon>0$,
	$$ \nu \left\{ x \in M \colon \left| \frac{1}{n} \, S_n \, \varphi (x) - \int_M \, \varphi \, d \nu  \right| > \epsilon \right\}   \to 0 \quad \text{ as } \ n \to \infty,$$
where $S_n \varphi := \varphi + \varphi \circ f + \cdots + \varphi \circ f^{n-1}$. Moreover, the precise convergence rate is quite important and usually only exponential and sub-exponential are useful for applications.

For a Markov system, we need a "stochastic" Birkhoff sum, and this is realized by Markov chain.

Let $X^+=M^\N$, consider the sequence of random variables $\{Z_n:X^+\to M\}_{n\in\N}$, $Z_n(x^+):= x_n$ where $x^+=\{x_n\}_{n\in \N}\in X^+$.

By Kolmogorov, given $\pi\in \Prob(M)$  there exists a unique
probability measure $\Pp_\pi$ on $X^+$ for which $\{Z_n\}_{n\in \N}$ is a  Markov chain with transition probability kernel $K$ and initial probability distribution $\pi$, i.e.,
such that for every $x\in M$, every Borel set $A\subset M$ and any $n\geq 1$,

\begin{enumerate}
	\item [(a)] $\Pp_\pi[ \, Z_n\in A \,\vert \,  Z_0,Z_1,\ldots, Z_{n-1}=x ] = K_{x}(A)$,
	\item [(b)] $\Pp_\pi[ \,Z_0 \in A]=\pi(A)$.
\end{enumerate}

When $\pi=\delta_x$, write $\Pp_x$. 

When $\pi=\nu$, $\Pp_\nu$ is shift invariant and makes
$\{Z_n\}_n$ a stationary process. It is called the Markov measure.

	$$ \Pp_\nu (B)=\int_M \Pp_x(B)\, d\nu(x)\quad \text{ and } \quad
	\EE_\nu[\psi]=\int_M \EE_x[ \psi ]\, d\nu(x)  $$
for any Borel set $B\subset X^+$ and any bounded measurable function
$\psi:X^+\to\R$.

Let $\{Z_n\}_{n\geq 0}$ be the $K$-Markov chain: $Z_n \colon X^+ \to M$, $Z_n(x^+)=x_n$. 
For an observable $\varphi \colon M \to \R$ and an index $j \ge 0$ let $$\varphi_j := \varphi \circ Z_j \colon X^+ \to \R.$$

Denote by 
$$
S_n \varphi := \varphi_0 + \cdots + \varphi_{n-1} =  \varphi(Z_0)+\cdots + \varphi(Z_{n-1})
$$
the corresponding ``stochastic'' Birkhoff sums. 

Here comes the main theorem of this section.

\begin{theorem}\label{absldt}
	Let $(M, K,\nu, \Escr)$ be a strongly mixing Markov system with mixing rate $r_n=\frac{1}{n^p}$, $p>0$. Then for all $\epsilon >0$ and $\varphi\in \Escr$ there are $c (\epsilon) > 0$ and $n (\epsilon) \in \N$ such that for all $n \ge n (\epsilon)$   we have
		$$
		\Pp_{\nu}\left\{\abs{\frac{1}{n}S_n\varphi-\int_M \varphi d\nu}>\epsilon\right\}\leq 8 e^{-c(\epsilon)n}
		$$
	where $c(\epsilon) = c_1 \, \epsilon^{2+\frac{1}{p}}$ and $n (\epsilon) = [c_2 \, \epsilon^{- \frac{1}{p}}]$ for constants $c_1 > 0$ and $c_2>0$.

If the strongly mixing condition is with the uniform norm $\norm{\cdot}_0$ on the left hand side, then the result holds with $\Pp_x$ for any $x\in M$.
\end{theorem}

The proof relies on the Bernstein's trick, the H\"older inequality and a trick of ``sparse'' rearrangement.

As an application, we proved an LDT (and also a CLT) for three levels of random torus translation realized by Markov chains for H\"older observables depending on the zero-th coordinate, the past coordinates and the full coordinates under a prevalent (mixing Diophantine, see Definition \ref{mixdc}) condition on the measure.

This theorem, containing the best existing exponential LDT result, is not only interesting itself, but also fundamental to our study because our mixed model can be embedded properly into a Markov system. Moreover, under certain assumptions we can prove that it is strongly mixing. This ensures the applicability of Theorem \ref{absldt} and the regularity of LE is thus in our hands thanks to the ACT.

\section{Regularity of the Lyapunov exponent}\label{example}
In this section, we state respectively two theorems for one-frequency and random-frequency mixed random-quasiperiodic cocycles which will appear in a forthcoming paper~\cite{CDK-paper4}.

Let $X:=\Gscr^\Z$ be the space of sequences
of  quasiperiodic cocycles and denote by
$\sigma\colon X\to X$ the usual  two-sided shift.

Recall that a compactly supported measure 
$\muh\in \Prob_c(\Gscr)$
determines a {\em random-quasiperiodic cocycle}
$F_\muh:X\times\T^d\times\R^m 
\to X\times\T^d\times\R^m $
by the formula
$$ F_\muh(\omega,\theta, v):= \left(\sigma\omega, \theta+\alpha, A(\theta)\, v \right),$$
where $\omega_0=(\alpha,A)$ is the  $0$-th coordinate  of the sequence $\omega\in X$.
The base  of this cocycle is the map
$f :X\times\T^d\to X\times\T^d$ defined by
$f(\omega, \theta):=(\sigma \omega, \theta+\alfa(\omega))$
where $\alfa:X\to\T^d$, $\alfa(\omega):=\alpha$, for $\omega_0=(\alpha,A)$ as above. The map $f$  preserves the measure $\muh^\Z\times m$, where $m$ is the  normalized Haar measure of $\T^d$.
The matrix valued function of this cocycle is
$\Ascr:X\times\T^d \to \SL_m(\R)$, defined by $\Ascr(\omega,\theta):=
A(\theta)$, where $\omega_0=(\alpha, A)$.

We say that the random-quasiperiodic cocycle determined
by $\muh \in \Prob_c(\Gscr)$ has a {\em single  frequency} $\alpha\in\T^d$ when
$\supp(\mu)\subset \Gscr_\alpha:=\{\alpha\}\times C^1(\T^d,\SL_m(\R))$. In this case the measure
$\muh \in \Prob_c(\Gscr_{\alpha})$ takes the form
$\muh=\delta_\alpha\times \nuh$ 
for some  $\nuh\in \Prob_c(C^1(\T^d,\SL_m(\R)))$.

We say that $\alpha\in\T^d$ is {\em rationally independent} if
$ \langle k, \alpha \rangle \notin \Z$
for all $k\in\Z^d\setminus\{0\}$.
Given  $L<\infty$  denote by $C^1_L(\T^d,\SL_m(\R))$   the 
set of all functions $A\in C^1(\T^d,\SL_m(\R)) $
such that $\norm{A}_{C^1} \leq L$ and $\norm{A^{-1}}_{C^1} \leq L $.

Finally let  $L_1(\muh)\geq L_2(\muh)\geq \ldots \geq L_d(\muh)$
denote the Lyapunov exponents of the random-quasiperiodic cocycle determined by
$\mu$. 

The main results of this section are

\begin{theorem}
	\label{thmA}
	Given $0<L<\infty$, $\alpha\in\T^d$ rationally independent
	and a measure $\nuh_0 \in \Prob_c(C^1_L(\T^d,\SL_m(\R)))$ assume that $\muh_0:=\delta_\alpha\times\nuh_0$ satisfies:
	\begin{enumerate}
		\item $\muh_0$ is quasi-irreducible;
		\item $L_1(\muh_0)> L_2(\muh_0)$.
	\end{enumerate}
	Then there exists $\delta>0$ such that the ball $\mathscr{U}$ of radius $\delta$  around $\nuh_0$  in
	the space $\Prob_c(C^1_L(\T^d,\SL_m(\R)))$, w.r.t. the Wasserstein distance, satisfies:
	\begin{enumerate}
		\item[(a)] The cocycles $F_\muh$ with  $\muh:=\delta_\alpha\times\nuh$ and $\nuh\in \mathscr{U}$, satisfy uniform LDT estimates of exponential type over $\mathscr{U}$.

		\item[(b)] The function $\mathscr{U}\in \nuh\mapsto L_1(\muh)$, where $\muh:=\delta_\alpha\times\nuh$,
		is H\"older continuous.
	\end{enumerate}
\end{theorem}

\bigskip

\begin{definition}\label{mixdc}
	We say that $\mu\in\Prob(\T^d)$ satisfies a mixing Diophantine condition if there are positive constants $\gamma$ and $\tau$
	such that
	\begin{equation}
		\label{mixing DC}
		\vert \hat \mu(k)\vert \leq  1-\frac{\gamma}{\norm{k}^\tau}\quad \forall k\in\Z^d\setminus\{0\}.
	\end{equation}
\end{definition}
Let $\rm{DC}_m(\gamma,\tau)$ denote the set of all probability measures $\mu\in\Prob(\T^d)$ satisfying~\eqref{mixing DC}.
Notice that $\rm{DC}_m(\gamma,\tau)$ is a compact and convex set of measures.

For any $L<\infty$ let $\Sigma_L:=\T^d\times C^1_L(\T^d,\SL_m(\R))$.

\begin{theorem}
	\label{thmB}
	Given positive constants $L$, $\gamma$ and $\tau$ 
	and a measure  $\muh_0 \in \Prob_c(\Sigma_L)$ satisfying:
	\begin{enumerate}
		\item $\muh_0$ is quasi-irreducible;
		\item $L_1(\muh_0)> L_2(\muh_0)$;
		\item $\mu_0:= \alfa_\ast \muh_0\in \rm{DC}_m(\gamma, \tau)$.
	\end{enumerate}
	Then there exists $\delta>0$ such that the ball $\mathscr{U}$ of radius $\delta$  around $\muh_0$  in
	the space $\Prob_c(\Sigma_L)\cap \alfa^{-1}(\rm{DC}_m(\gamma,\tau))$, w.r.t. the Wasserstein distance, satisfies:
	\begin{enumerate}
		\item[(a)] The cocycles $F_\muh$ with $\muh\in \mathscr{U}$, satisfy uniform LDT estimates of exponential type over $\mathscr{U}$.
		
		\item[(b)] If $\tau$ is small enough, all cocycles $F_\muh$ with $\muh\in \mathscr{U}$, satisfy a CLT.	
		
		\item[(c)] The function $\mathscr{U}\in \muh\mapsto L_1(\muh)$ 	is H\"older continuous.
	\end{enumerate}
\end{theorem}

As explained previously, Theorem \ref{thmA} and \ref{thmB} will be consequences of the abstract LDT theorem. However, there is an essential difference between one frequency and random frequencies. That is, in the one-frequency case, we do not have the ``full'' strong mixing of our Markov operator as the torus translation is never mixing. Luckily, we may separate the eigenspace into a direct sum of ``quasi-periodic'' subspace and a ``random'' one where the action of our Markov operator becomes a torus translation on the quasi-periodic subspace which has no loss of parameters (due to uniform convergence of Birkhoff sums by unique ergodicity) and thus LDT is immediate. For the random part, we can prove strong mixing so the LDT holds also. Overall, although the ``full'' strong mixing does not hold,  the ``full'' LDT does hold.

The technical part after embedding our model into some suitable Markov system lies in proving the direct sum decomposition and the strong mixing of the Markov operator restricted to the  random part (which does not mean that the random-frequency case is easier). In fact, lots of effort was devoted to both cases,  involving many results from the previous papers \cite{CDK-paper2,CDK-paper1,CDK-paper3}.

Definitely, Theorem \ref{thmA} and \ref{thmB} apply to the Schr\"odinger operators in Propositions \ref{intro prop1} \ref{intro prop2}\ref{intro prop3}, which again shows that randomness dominates quasi-periodicity in the sense of regularity of the Lyapunov exponent. For example, the Lyapunov exponent of a quasi-periodic cocycle can be discontinuous even with a very good topology~\cite{WangYou}. However, as the randomness participates, the mixed cocycle exhibits H\"older continuous Lyapunov exponent consistent with the random case.

To finish this section, we note that H\"older continuity with respect to the measure is the strongest result as it implies the same result for the energy, the frequency and the potential about which many mathematicians in Mathematical Physics are more concerned.

\section{Stability of the Lyapunov exponent}\label{uniform}
With all the preparations in hand, we are finally ready to hit the original question of You. As we may have seen, it takes us a very long period as well as lots of effort to stand finally at the equal stage of You's question. And still, this is just the beginning of solving it.

\subsection*{Formulation of the problem}

Given a rationally independent frequency $\alpha$ on the $d$-dimensional torus $\T^d = (\R/\Z)^d$, let $\tau_\alpha (\theta) = \theta + \alpha \, \mod 1$ be the translation by $\alpha$ on $\T^d$ endowed with the Haar measure $m$.

Given a function $v \colon \T^d \to \R$, let $A_{0, v} \colon \T^d \to \SL_2 (\R)$,
$$ A_{0, v} (\theta):= \begin{bmatrix}
	v (\theta) & -1  \\ 1 & 0
\end{bmatrix}$$

The skew-product map $F_{0, v} \colon \T^d \times \R^2 \to \T^d \times \R^2$,
$$F_{0, v} (\theta, p) := (\tau_\alpha (\theta), A_0 (\theta) p )$$
is called a quasiperiodic Schr\"odinger cocycle. 
We will also refer to the matrix valued function $A_0$ as being the same quasiperiodic cocycle.

Let $\rho \in \Prob_c (\R)$ and consider the following random perturbation of the quasiperiodic cocycle $A_0$:
$$A_{\ep, v} (\om_0, \theta) := 
\begin{bmatrix}
	v (\theta) + \ep \, \om_0  & -1  \\ 1 & 0
\end{bmatrix} = P (\ep \, \om_0) \, A_{0, v} (\theta) ,
$$
where
$$ P (\om_0):= \begin{bmatrix}
	1 & \om_0  \\ 0 & 1
\end{bmatrix}$$

Denote by $\Sigma$ the support of $\rho$, let $X := \Sigma^\Z$ and consider the two-sided Bernoulli shift $\sigma$ on $(X, \rho^\Z)$. 
The random perturbation $A_{\ep, v}$ determines the linear cocycle $F_{\ep, v} \colon X \times \T^d \times \R^2 \to X \times \T^d \times \R^2$,
$$ F_{\ep, v} (\om, \theta, p) := (\sigma \om, \theta + \alpha, A_{\ep, } (\om, \theta) p )$$
over the mixed random-quasiperiodic base dynamics 
$$X \times \T^d \times \ni (\om, \theta) \mapsto f (\om, \theta) = (\sigma \om, \tau_\alpha (\theta)) \, .$$

For this Schr\"odinger case, we have the following theorem.

\begin{theorem}\label{main thm Schrodinger}
	Assume  the function $v \colon \T^d \to \R$ is analytic,  the frequency $\alpha \in \T^d$ satisfies a Diophantine condition,
	$L_1(A_{0,v})>0$ and   $\rho \in \Prob_c (\R)$ has positive pointwise Hausdorff dimension, i.e.,  for some constants $C< \infty$ and $\vartheta > 0$
	$$\rho (B (x, r)) \le C \, r^\vartheta\quad  \forall x \in \R, r > 0 .$$ 	
	
	Then the map
	$$ (\ep, v) \mapsto L_1 (F_{\ep, v}) $$
	is locally weak-H\"older continuous on $[0, \infty) \times C_r^\om (\T^d, \R)$. 
	In particular the Lyapunov exponent is stable under random perturbations. 
\end{theorem} 

This theorem is a particular version of the following general result.

\newcommand{\slin}{\mathrm{sl}}

\subsection*{Definitions and general statement}

Let $\slin_m(\R)$ denote the Lie algebra of the group
$\SL_m(\R)$ and consider a setting with the following four ingredients:

\begin{enumerate}
	\item  $\alpha\in\T^d$ a frequency satisfying a Diophantine condition;
	\item $A_0\in C^\omega_r(\T^d, \SL_m(\R))$ an analytic cocycle;
	\item $\mu\in\Prob_c(\slin_m(\R))$ a compactly supported measure.
\end{enumerate}

These determine the following objects:
A Bernoulli space of sequences $\left( X:=\slin_m(\R)^\Z, \mu^\Z\right)$ endowed with a two-sided Bernoulli shift $\sigma:X\to X$. With it we define the base map $f:X\times\T^d\to X\times \T^d$,
$f(\omega,\theta):=(\sigma\omega, \theta+\alpha)$,
which together with the  matrix valued function
$A_\epsilon:X\times\T^d\to \SL_m(\R)$, $A_\epsilon(\omega,\theta):= e^{\epsilon\, \omega_0}  A_0(\theta)$ determines the randomly perturbed linear cocycle	 
$F_{\alpha, A_0, \mu, \epsilon}:X\times \T^d\times\R^m\to X\times \T^d\times\R^m$, 
$$F_{\alpha, A_0, \mu, \epsilon}(\omega, \theta, p):=(f(\omega, \theta),  A_\epsilon(\omega,\theta)\, p) . $$

We make two assumptions. First we assume a gap on the Lyapunov spectrum.

\begin{enumerate}
	\item[(H1)]  $L_1(\alpha, A_0)>L_2(\alpha, A_0)$.
\end{enumerate}

The second hypothesis is a bit more technical.
\begin{enumerate}
	\item[(H2)]  There exist constants 
	$C<\infty$, $\vartheta>m-2$, $\varepsilon_0>0$ and
	$k_0\in\N$ such that for every $\theta\in\T^d$,
	$\hat p, \hat q\in\Proj$ and $0<r<\varepsilon<\varepsilon_0$ 
	$$ \mu^\Z\left\{ \omega\in X\colon \, \hat A_{\varepsilon_0}^{k_0}(\omega,\theta)\,\hat p \in B(\hat q, r)\, \right\} < C\, \left(\frac{r}{\varepsilon} \right)^\vartheta . $$
\end{enumerate}

Then we have a more general theorem as follows.

\begin{theorem}
	Fix $\alpha\in\T^d$ satisfying a Diophantine condition and  the space $\mathscr{U}\subseteq \Prob(\slin_m(\R))\times C^\omega_r(\T^d, \SL_m(\R))\times [0, +\infty)$ of all tuples $(\mu, A_0, \epsilon)$ such that
	$(\mu, A_0)$ satisfy hypothesis (H1) and (H2).
	Then the function $L_1:\mathscr{U}\to\R$, $(\mu,A_0,\epsilon)\mapsto L_1(F_{\alpha, A_0, \mu, \epsilon})$
	is locally weak-H\"older continuous.
\end{theorem}

The proof is based on an approximation lemma by quasi-periodic estimates, results of mixed random-quasiperiodic cocycles and a bridging argument by the Avalanche Principle (AP). The key idea is that when there are not many iterates, the system is behaving more like the quasi-periodic one while if there are enough iterates, the randomness starts to work and dominate. Finally, by AP the hole in the middle can be perfectly filled in!

We believe that (H2) could be removed though at this point we can not. That is why we call our result partial. In fact, there are many other questions that we are trying to resolve. For example, stability in the random-frequency case and metal-insulator transition from quasi-periodic models to mixed models, etc.

Hopefully, the papers in preparation will be released soon.

\bigskip
Last but not least:

\subsection*{Acknowledgments} The author would like to give his deepest and sincerest gratitude to Pedro Duarte, Silvius Klein, Jiangong You and Qi Zhou for their consistent support and persistent inspiration.

\bibliographystyle{amsplain} 
\bibliography{references} 

\end{document}